\documentclass[fleqn]{article}
 \usepackage{amscd}
                     \usepackage{amssymb}
 \newcommand{\be}{\begin{equation}}
       \newcommand{\ee}{\end{equation}}
       \newcommand{\ba}{\begin{eqnarray}}
        \newcommand{\ea}{\end{eqnarray}}
 \newcommand{\ban}{\begin{eqnarray*}}
 \newcommand{\ean}{\end{eqnarray*}}

 \newcommand{\lp}{\langle}
 \newcommand{\rp}{\rangle}
 \newcommand{\ra}{\rightarrow}

 \newcommand{\sect}[1]{\section{#1} \setcounter{equation}{0}}

 \newcommand{\vol}{\mathrm{Vol}}
 \newcommand{\diam}{\mathrm{diam}}
 \newcommand{\Ric}{\mathrm{Ric}}
 \newcommand{\Hess}{\mathrm{Hess}}
  \newcommand{\tr}{\mbox{tr}}

 \newtheorem{theo}{Theorem}[section]

\begin{document}
 \newtheorem{defn}[theo]{Definition}
 \newtheorem{ques}[theo]{Question}
 \newtheorem{lem}[theo]{Lemma}
 \newtheorem{lemma}[theo]{Lemma}
 \newtheorem{prop}[theo]{Proposition}
 \newtheorem{coro}[theo]{Corollary}
 \newtheorem{ex}[theo]{Example}
 \newtheorem{note}[theo]{Note}
 \newtheorem{conj}[theo]{Conjecture}

\title{Manifolds with A Lower Ricci Curvature Bound}
\author{Guofang Wei \thanks{Partially supported by NSF grant DMS-0505733}}
\date{}
\maketitle
\begin{abstract}
This paper  is a survey on the structure of manifolds with a lower Ricci curvature bound. 
\end{abstract}

\sect{Introduction} The purpose of this paper is to give a survey on
the structure of manifolds with a lower Ricci curvature bound. A Ricci
curvature bound is weaker than a sectional curvature bound but stronger
than a scalar curvature bound. Ricci curvature is also special that it occurs in the Einstein
equation and in the Ricci flow. The study of manifolds
with lower Ricci curvature bound has experienced tremendous progress
in the past fifteen years. Our focus in this  article is strictly restricted to results with only
Ricci curvature bound, and no result with sectional curvature bound
is presented unless for straight comparison. The reader is referred
to John Lott's article in this volume for the recent important
development concerning Ricci curvature for metric measure spaces by Lott-Villani and Sturm.  We start
by introducing the basic tools for studying manifolds with lower
Ricci curvature bound (Sections 2-4), then discuss the structures of
these manifolds (Sections 5-9), with examples in Section 10.

The most basic tool in studying manifolds with Ricci curvature bound
is the Bochner formula. From there one can derive  powerful
comparison tools like the mean curvature comparison, the Laplacian
comparison, and the relative volume comparison.  For the Laplacian
comparison (Section 3) we discuss the global version in three weak
senses (barrier, distribution, viscosity) and clarify their
relationships (I am very grateful to my colleague Mike Crandall for many
helpful discussions and references on this issue).
 A generalization of the volume comparison theorem  to an integral Ricci curvature bound is also presented (Section 4).  
 Important tools such as Cheng-Yau's gradient estimate and Cheeger-Colding's segment inequality are presented in Sections 2
  and 4 respectively. Cheeger-Gromoll's splitting theorem and Abresch-Gromoll's excess estimate  are presented in Sections 5 and 8 respectively.

From comparison theorems, various quantities like the volume, the
diameter, the first Betti number, and the first eigenvalue are
bounded by the corresponding quantity of the model. When equality occurs one has the
rigid case. In Section 5 we discuss many rigidity and stability
results for nonnegative and positive Ricci curvature. The Ricci
curvature lower bound gives very good control on the fundamental
group and the first Betti number of the manifold; this is covered in
Section 6 (see also the very recent survey article by Shen-Sormani
\cite{Shen-Sormani-survey} for more elaborate discussion). In
Sections~\ref{Gromov-Hausdorff}, 8, and 9 we discuss rigidity and
stability for manifolds with lower Ricci curvature bound under
Gromov-Hausdorff convergence, almost rigidity results, and the
structure of the limit spaces, mostly due to Cheeger and Colding.
Examples of manifolds with positive Ricci curvature are presented in
Section 10.

Many of the results in this article are covered in the very nice
survey articles \cite{Zhu1997, Cheeger2001}, where complete proofs
are presented. We benefit greatly from  these two articles. Some materials here are adapted directly from \cite{Cheeger2001} and we are very grateful to Jeff Cheeger for his permission. We also benefit from
\cite{Gromoll, Cheeger2002} and  the lecture notes
\cite{Wei-lecture-notes} of  a topics course I taught at UCSB. I
would also like to thank Jeff Cheeger, Xianzhe Dai, Karsten Grove, Peter Petersen,
Christina Sormani, and William Wylie for reading earlier versions of this article and for their
helpful suggestions.

  \sect{Bochner's formula and the mean
curvature comparison}

For a smooth function $f$ on a complete Riemannian manifold $(M^n,g)$, the gradient of f is the vector field $\nabla f$ such that $\lp \nabla f, X \rp = X(f)$ for all vector fields $X$ on $M$. The Hessian of $f$ is the symmetric bilinear form
\[
\Hess (f) (X,Y) = XY (f) - \nabla_XY(f) = \lp \nabla_X \nabla f,Y \rp,
\]
and the Laplacian is the trace $\Delta f = \tr (\Hess f)$. For a bilinear form $A$, we denote $|A|^2 = \tr (AA^t)$.  The most basic tool in studying manifolds with Ricci curvature bound is the Bochner formula. Here we state the formula for functions.
\begin{theo}[Bochner's Formula] For a smooth function $f$ on a complete Riemannian manifold $(M^n,g)$,
\be  \label{bochner}
\frac 12 \Delta |\nabla f|^2 = |\Hess f|^2 + \lp \nabla f, \nabla (\Delta f) \rp + \Ric (\nabla f, \nabla f).
\ee
\end{theo}

This formula has many applications. In particular, we can apply it to the distance function, harmonic functions, and the eigenfunctions among others.  The formula has a more general version (Weitzenb\"{o}ck type) for vector fields (1-forms), which also works nicely on Riemannian manifolds with a smooth measure \cite{Lott2003, Perelman-math.DG/0211159} where Ricci and all adjoint operators are defined with respect to the measure.

Let $r(x) = d(p,x)$ be the distance function from $p \in M$. $r(x)$ is a Lipschitz function and is smooth on $M\setminus \{p, C_p\}$, where $C_p$ is the cut locus of $p$. At smooth points of $r$,
\be
|\nabla r| \equiv 1, \ \
\Hess \, r = II, \ \ \Delta r = m,
\ee
where $II$ and $m$ are the second fundamental form and mean curvature of the geodesics sphere $\partial B(p,r)$.

Putting $f(x) = r(x)$ in (\ref{bochner}), we obtain the Riccati equation along a radial geodesic,
\be
0 = | II |^2 + m'+ \Ric (\nabla r, \nabla r).
\ee
By the Schwarz inequality,
\[
| II |^2 \ge \frac{m^2}{n-1}.
\]
Thus, if $\Ric_{M^n} \ge (n-1)H$, we have the Riccati inequality,
\be  \label{Riccati-inequ}
m' \le -\frac{m^2}{n-1} - (n-1)H.
\ee

Let $M^n_H$ denote the complete simply connected space of constant curvature $H$ and $m_H$ the mean curvature of its geodesics sphere, then
\be \label{Riccati-inequ-H}
m'_H  =  -\frac{m_H^2}{n-1} - (n-1)H.
\ee

Since $\lim_{r \ra 0} (m - m_H) = 0$, using (\ref{Riccati-inequ}), (\ref{Riccati-inequ-H}) and the standard Riccati equation comparison, we have
\begin{theo}[Mean Curvature Comparison] If $\Ric_{M^n} \ge (n-1)H$, then along any minimal geodesic segment from $p$,
\be  \label{mean-comp}
m(r) \le m_H(r).
\ee
Moreover, equality holds if and only if all radial sectional
 curvatures are equal to $H$.
\end{theo}

By applying the Bochner formula to $f =  \log u$ with an appropriate cut-off function and looking at the maximum point  one has Cheng-Yau's gradient estimate for harmonic functions \cite{Cheng-Yau1975}.
\begin{theo}[Gradient Estimate, Cheng-Yau 1975]
Let $\Ric_{M^n} \ge (n-1)H$ on $B(p,R_2)$ and $u: B(p,R_2) \ra \mathbb R$ satisfying $u >0, \Delta u = 0$. Then for $R_1 < R_2$, on $B(p, R_1)$,
\be  \label{gradient}
\frac{|\nabla u|^2}{u^2} \le c(n,H,R_1,R_2).
\ee
\end{theo}
If $\Delta u = K(u)$, the same proof extends and one has \cite{Cheeger2001}
\be   \label{gradient2}
\frac{|\nabla u|^2}{u^2} \le \max \{ 2u^{-1} K(u), c(n,H,R_1,R_2)+ 2u^{-1}K(u) - 2K'(u)\}.
\ee

\sect{Laplacian comparison}

Recall that $m = \Delta r$. From (\ref{mean-comp}), we get the local Laplacian comparison for distance functions
\be  \label{lap-com}
\Delta r \le \Delta_H r,  \ \ \mbox{for all}\ x \in M\setminus \{p, C_p\}.
\ee
 Note that if $x \in C_p$, then either $x$ is a (first) conjugate point of $p$ or there are two distinct minimal geodesics connecting $p$ and $x$ \cite{Cheeger-Ebin}, so $x \in $\{conjugate locus of $p\} \cup \{$the set where $r$ is not differentiable\}. The conjugate locus of p consists of the critical values of exp$_p$. Since exp$_p$ is smooth, by Sard's theorem, the conjugate locus has measure zero. The set where $r$ is not differentiable has measure zero since $r$ is Lipschitz. Therefore the cut locus $C_p$ has measure zero. One can show  $C_p$ has measure zero more directly by observing that the region inside the cut locus is star-shaped \cite[Page 112]{Chavel1993}. The above argument has the advantage that it can be extended easily to show that Perelman's $l$-cut locus \cite{Perelman-math.DG/0211159}  has measure zero since the ${\mathcal L}$-exponential map is smooth and the $l$-distance function is locally Lipschitz.

In fact the Laplacian comparison (\ref{lap-com}) holds globally in various weak senses. First we review the definitions (for simplicity we only do so for the Laplacian) and study the relationship between these different weak senses.

For a continuous function $f$ on $M, q \in M$, a function $f_q$ defined in a neighborhood $U$ of $q$, is an upper barrier of $f$ at $q$ if $f_{q}$ is $C^2(U)$    and
 \be  \label{up-barrier}
f_{q} (q) = f(q), \ \ \
f_{q} (x) \ge f(x)  \  (x \in U).
\ee

\begin{defn}
For a continuous function $f$ on $M$, we say $\Delta f (q) \le  c$ in the barrier sense ($f$ is a barrier subsolution to the equation $\Delta f =c$ at $q$), if for all $\epsilon > 0$, there exists an upper barrier $f_{q, \epsilon} $ such that $\Delta f_{q, \epsilon} (q) \le c + \epsilon.$
\end{defn}
This notion was defined by Calabi  \cite{Calabi1958} back in 1958 (he used the terminology  ``weak sense" rather than  ``barrier sense").  
A weaker version is in the sense of viscosity,  introduced by Crandall and Lions in \cite{Crandall-Lions1983}.
\begin{defn}
For a continuous function $f$ on $M$, we say $\Delta f (q) \le  c$ in the viscosity sense ($f$ is a viscosity subsolution of $\Delta f =c$ at $q$), if $\Delta \phi (q) \le c$
whenever $\phi \in C^2(U)$ and $(f-\phi)(q) = \inf_{U} (f-\phi)$, where $U$ is a neighborhood of $q$.
\end{defn}

Clearly  barrier subsolutions are viscosity  subsolutions.

Another very useful notion is subsolution in the sense of distributions.
\begin{defn}
For  continuous functions $f, h$ on an open domain $\Omega \subset M$, we say $\Delta f  \le  h$ in the distribution sense ($f$ is a distribution subsolution of $\Delta f  =  h$) on $\Omega$, if
 $ \int_{\Omega} f \Delta \phi  \le \int_{\Omega} h  \phi$
for all $\phi \ge 0$ in $C^\infty_0(\Omega)$.
\end{defn}

By \cite{Ishii} if $f$ is a viscosity subsolution of $\Delta f  = h$
on $\Omega$, then it is also a distribution subsolution  and vice
verse, see also \cite{Lions}, \cite[Theorem 3.2.11]{Hormander1994}.

For geometric applications, the barrier and distribution sense are very useful and the barrier sense is often easy to check. Viscosity gives a bridge between them.  As observed by Calabi \cite{Calabi1958} one can easily construct upper barriers for the distance function.
\begin{lem}  If $\gamma$ is minimal from $p$ to $q$, then for all $\epsilon > 0$, the function $r_{q,\epsilon}(x) = \epsilon + d(x,\gamma(\epsilon))$, is an upper barrier for the distance function $r(x) = d(p,x)$ at $q$.
\end{lem}

Since $r_{q,\epsilon}$  trivially satisfies (\ref{up-barrier})
the lemma follows by observing that it is smooth in a neighborhood of $q$. 

Upper barriers for Perelman's $l$-distance function can be constructed very similarly. 

Therefore the Laplacian comparison (\ref{lap-com}) holds globally in
all the weak senses above.  Cheeger-Gromoll (unaware of Calabi's work at the time) had proved the Laplacian comparison in the distribution sense directly by observing  the very useful fact that near the cut locus $\nabla r$  points towards  the cut locus \cite{Cheeger-Gromoll1971}, see also \cite{Cheeger2001}.  (However it is not clear if this fact holds for Perelman's $l$-distance function.)

One reason why these weak subsolutions are so useful is that they
still satisfy the following classical Hopf strong maximum principle, see \cite{Calabi1958},  also e.g. \cite{Cheeger2001} for the barrier sense, see \cite{Littman1959, Kawohl-Kutev1998} for the distribution and viscosity  senses, also \cite[Theorem 3.2.11]{Hormander1994} in the Euclidean case.
\begin{theo}[Strong Maximum Principle] \label{maximum-p}
If on a connected open set, $\Omega \subset M^n$, the function $f$
has an interior minimum and  $\Delta f \le 0$ in any of the weak senses above, then
$f$ is constant on $\Omega$.
\end{theo}

 These weak solutions also enjoy the regularity
(e.g. if $f$ is a weak sub and sup solution of $\Delta f =0$, then $f$ is
smooth), see e.g. \cite{Gilbarg-Trudinger2001}.

The Laplacian comparison also works  for radial functions (functions composed with the distance function). In geodesic polar coordinate, we have
\be
\Delta f  =  \tilde{\Delta} f + m(r, \theta) \frac{\partial}{\partial r}f +
 \frac{\partial^2 f}{\partial r^2},
\ee
where $\tilde{\Delta}$ is the induced Laplacian on the sphere and
$m(r,\theta)$ is the mean curvature of the geodesic sphere in the inner normal
direction. Therefore
\begin{theo}[Global Laplacian Comparison] \label{laplacian-comp}
If $\Ric_{M^n} \ge (n-1)H$, in all the weak senses above, we have
\ba
\Delta f(r) & \le & \Delta_H f(r)   \ \ \  (\mbox{if}\ f' \ge 0), \\
\Delta f(r) & \ge & \Delta_H f(r)   \ \ \  (\mbox{if}\ f' \le 0).
\ea
\end{theo}

\sect{Volume comparison}
For $p \in M^n$, use exponential polar coordinate around $p$ and write the volume element $d\, vol = \mathcal A(r,\theta) dr \wedge d\theta_{n-1}$, where $d\theta_{n-1}$ is the standard volume element on the unit sphere $S^{n-1}(1)$. By the first variation of the area (see \cite{Zhu1997})
\be \label{area}
\frac{\mathcal A'}{\mathcal A} (r,\theta) = m (r, \theta).
\ee
Similarly, define ${\mathcal A}_H$ for the model space $M_H^n$. The mean curvature comparison and  (\ref{area}) gives
the volume element comparison. Namely if $M^n$ has $\Ric_M \geq (n-1)H$, then
\be  \label{vol-elem}
\frac{\mathcal{A}(r,\theta)}{\mathcal{A}_H
  (r,\theta)} \ \mbox{ is nonincreasing along any minimal geodesic segment from} \ p.
  \ee

Integrating (\ref{vol-elem}) along the sphere directions, and then the radial direction gives the relative area and volume comparison, see e.g. \cite{Zhu1997}.
\begin{theo}[Bishop-Gromov's Relative Volume Comparison] \label{Bishop-Gromov}
Suppose $M^n$ has $\Ric_M \geq (n-1)H$.  Then
\be  \label{area-vol}
\frac{\vol \, (\partial B(p,r))}{\vol_H (\partial B(r))} \ \mbox{and}\  \frac{\vol\, (B(p,r))}{\vol_H (B(r))}  \ \mbox{ are nonincreasing in}  \ r.
\ee
In particular,
\be
\vol \, (B(p,r)) \le \vol_H (B(r)) \ \ \  \mbox{for all } \ r > 0,  \label{vol-absolute}
\ee
\be
 \frac{\vol\, (B(p,r))}{\vol\, (B(p, R))}   \ge  \frac{\vol_H (B(r))}{\vol_H(B(R))} \ \ \ \mbox{for all } \ 0 < r \le R,
 \label{relative-vol} \ee
  and equality holds if and only if $B(p,r)$ is isometric to $B_H(r)$.
\end{theo}

 This is a powerful result because  it is a global comparison. The volume of any ball is bounded above by the volume of
 the corresponding ball in the model, and if  the volume of a big ball has a lower bound, then all smaller balls also
 have lower bounds. One can also apply it to an annulus or a section of the
 directions. For topological applications see Section 6.

The volume element comparison (\ref{vol-elem}) can also be used to prove a heat kernel comparison \cite{Cheeger-Yau1981} and Cheeger-Colding's segment inequality \cite[Theorem 2.11]{Cheeger-Colding1996}, see also \cite{Cheeger2001}.

Given a function $g \ge 0$ on $M^n$, put
\[
\mathcal F_g(x_1,x_2) = \inf_\gamma \int_0^l  g(\gamma(s)) ds,
\]
where the inf is taken over all minimal geodesics $\gamma$ from $x_1$ to $x_2$ and $s$ denotes the arclength.
\begin{theo}[Segment Inequality, Cheeger-Colding 1996] \label{thm-segment} Let $\Ric_{M^n} \ge -(n-1)$, $A_1, A_2 \subset B(p,r)$, and $r \le R$. Then
\be  \label{segment}
\int_{A_1\times A_2} \mathcal F_g(x_1,x_2) \le c(n,R) \cdot  r \cdot  (\vol (A_1) + \vol (A_2)) \cdot \int_{B(p,2R)} g,
\ee
where $c(n,R) = 2 \sup_{0<\frac s2 \le u \le s, 0<s<R} \frac{\vol_{-1} (\partial B(s))}{\vol_{-1} (\partial B(u))}$.
\end{theo}

The segment inequality shows that if the integral of $g$ on a ball is small then the integral of $g$ along almost all segments is small. It also implies a Poincar\'e inequality of type $(1,p)$ for all $p \ge 1$ for manifolds with lower Ricci curvature bound  \cite{Buser1982}. In particular it gives a lower bound on the first eigenvalue of the Laplacian for the Dirichlet problem on a metric ball; compare \cite{Li1993}.

 The volume comparison theorem can be generalized to an integral Ricci lower bound \cite{Petersen-Wei1997}, see also \cite{Gallot1988, Yang1992}. For convenience we introduce some notation.

 For each $x\in M^n$ let $\lambda\left( x\right) $ denote the smallest
eigenvalue for the Ricci tensor $\mathrm{Ric}:T_{x}M\rightarrow
T_{x}M,$ and $\mathrm{Ric}_-^H(x) = \left( (n-1)H - \lambda
(x)\right)_+ = \max \left\{ 0, (n-1)H - \lambda (x) \right\}$. Let
 \be \| \mathrm{Ric}_-^H \|_p (R) = \sup_{x\in M} \left( \int_{B\left( x,R\right) } (\mathrm{Ric}_-^H)^{p}\,
  d vol\right)^{\frac 1p}. \ee
 $\| \mathrm{Ric}_-^H \|_p $ measures the amount of Ricci
curvature lying below $(n-1)H$ in the $L^p$ sense. Clearly $\|
\mathrm{Ric}_-^H \|_p (R) = 0$ iff $\Ric_M \ge (n-1)H$.

Parallel to the mean curvature comparison theorem (\ref{mean-comp}) under
pointwise Ricci curvature lower bound, Petersen-Wei
\cite{Petersen-Wei1997} showed one can estimate the amount of mean
curvature bigger than the mean curvature in the model by the amount of Ricci
curvature lying below $(n-1)H$ in $L^p$ sense. Namely for any $p > \frac
n2$, $H \in \mathbb R$, and  when
$H>0$ assume $r\le \frac{\pi}{2\sqrt{H}}$, we have \be  \label{mean-estimate} \left( \int_{B\left(
x,r\right) }\left( m-m_H \right) _{+}^{2p} \,d vol
\right)^{\frac{1}{2p}} \leq C\left( n,p\right) \cdot \left( \|
\mathrm{Ric}_-^H \|_p (r) \right)^{\frac 12}. \ee

Using (\ref{mean-estimate}) we have
\begin{theo}[Relative Volume Estimate, Petersen-Wei 1997]  Let $x \in M^n, H \in \mathbb R$ and $p>\frac n2$ be given, then there is a constant $C(n,p,H,R)$ which is nondecreasing in $R$ such that if $r\leq R$ and when $H>0$ assume that $R \le \frac{\pi}{2\sqrt{H}}$  we have
\be  \label{vol-estimate} \left( \frac{\vol \,B\left( x,R\right)
}{\vol_H\left(B(R)\right) }\right) ^{\frac 1{2p}}- \left( \frac{\vol
\,B\left( x,r\right) }{\vol_H  \left(B(r)\right) }\right) ^{\frac
1{2p}} \leq C\left( n,p,H ,R\right) \cdot \left( \| \mathrm{Ric}_-^H
\|_p (R) \right)^{\frac 12}. \ee
 Furthermore when $r=0$ we obtain
 \be
\vol \,B\left( x,R\right) \leq \left( 1+ C\left( n,p,H,R\right)
\cdot \left( \| \mathrm{Ric}_-^H \|_p (R) \right)^{\frac 12}
\right)^{2p} \vol_H\left( B(R)\right). \ee
\end{theo}

Note that when $\| \mathrm{Ric}_-^H \|_p (R) =0$, this gives the
Bishop-Gromov relative volume comparison.

Volume comparison is a powerful tool for studying manifolds with
lower Ricci curvature bound  and has many applications. As a result
of (\ref{vol-estimate}), many results with pointwise Ricci lower
bound (i.e. $\| \mathrm{Ric}_-^H \|_p (r)=0$) can be extended to the
case when $\| \mathrm{Ric}_-^H \|_p (r)$ is very small
\cite{Gallot1988,Petersen-Wei1997, Petersen-Sprouse1998,
Dai-Petersen-Wei2000, Sprouse2000, Petersen-Wei2001, Dai-Wei2004,
Aubry-preprint}.

Perelman's reduced volume monotonicity \cite{Perelman-math.DG/0211159}, a basic and
powerful tool in his work on Thurston's geometrization conjecture,
is a generalization of Bishop-Gromov's volume comparison to Ricci
flow. In fact Perelman gave a heuristic argument that volume
comparison on an infinite dimensional space (incorporating the Ricci
flow) gives the reduced volume monotonicity. It would be very interesting
to investigate this relationship further.

\sect{Rigidity results and stability}

From comparison theorems, various quantities are bounded by that of the
model. When equality occurs one has the rigid case. In this section we concentrate on the rigidity and stability results for nonnegative and positive Ricci curvature. See Section~\ref{Gromov-Hausdorff} for rigidity and  stability under Gromov-Hausdorff convergence and a general lower bound.

The simplest rigidity is the maximal volume. From the equality of volume comparison (\ref{vol-absolute}), we deduce that if $M^n$ has $\Ric_M \ge
n-1$ and $\vol_M = \vol (S^n)$, then $M^n$ is isometric to $S^n$. Similarly if $M^n$ has $\Ric_M \ge 0$ and $\lim_{r \to \infty}  \frac{\vol B(p,r)}{\omega_n r^n} =1$, where $p \in M$ and $\omega_n$ is the volume of the unit ball in  $\mathbb R^n$, then $M^n$ is isometric to $\mathbb R^n$.

From the equality of the area of geodesic ball (the first quantity in (\ref{area-vol})) we get another volume rigidity: volume annulus implies metric annulus. This is first observed in \cite[Section 4]{Cheeger-Colding1996}, see also \cite[Theorem 2.6]{Cheeger2002}. For the case of nonnegative Ricci curvature, this result says that if $\Ric_{M^n} \ge 0$ on the annulus $A(p, r_1, r_2)$, and
\[
\frac{\vol (\partial B(p, r_1))}{\vol (\partial B(p, r_2))} = \frac{r_1^{n-1}}{{r_2}^{n-1}},
\]
 then the metric on $A(p, r_1, r_2)$ is of the form, $dr^2 + r^2 \tilde{g}$, for some smooth Riemannian metric $\tilde{g}$ on $\partial B(p, r_1)$.

By Myers' theorem (see Theorem~\ref{Myers}) when Ricci curvature has a positive lower bound the diameter is bounded by the diameter of the model.  In the maximal case, using an eigenvalue comparison (see below) Cheng \cite{Cheng1975} proved that if $M^n$ has $\Ric_M \ge n-1$ and $\diam_M = \pi$, then $M^n$ is isometric to $S^n$.
This result can also be
directly proven using volume comparison \cite{Shiohama1983,Zhu1997}.

The maximal diameter theorem for the noncompact case is given by Cheeger-Gromoll's splitting
theorem \cite{Cheeger-Gromoll1971}. The splitting theorem is the most important rigidity result, it plays a very important
role in studying manifolds with nonnegative Ricci curvature and
manifolds with general Ricci lower bound.
\begin{theo}[Splitting Theorem, Cheeger-Gromoll 1971]  \label{splitting} Let $M^n$ be
a complete Riemannian manifold with $\Ric_M \ge 0$. If $M$ has a
line, then $M$ is isometric to the product $\mathbb R \times
N^{n-1}$, where $N$ is an $n-1$ dimensional manifold with $\Ric_N
\ge 0$.
\end{theo}

The result can be proven using the global Laplacian comparison
(Theorem~\ref{laplacian-comp}), the strong maximum principle
(Theorem~\ref{maximum-p}), the Bochner formula (\ref{bochner}) and the de
Rham decomposition theorem, see e.g. \cite{Zhu1997, Cheeger2001,
Petersen-book} for detail.

As an application of the splitting theorem we have that the first
Betti number of $M$ is less than or equal to  $n$ for $M^n$ with $\Ric_M \ge
0$, and  $b_1 =n$ if and only if $M$ is isometric to $T^n$ (the flat
torus).

Applying the Bochner formula  (\ref{bochner}) to the first eigenfunction Lichnerowicz showed that if $M^n$ has $\Ric_M \ge n-1$, then
the first eigenvalue $\lambda_1(M) \ge n$ \cite{Lichnerowicz1958}. Obata showed that if  $\lambda_1(M) =n$ then $M^n$ is isometric to $S^n$ \cite{Obata1962}.

From these rigidity results (the equal case), we naturally ask what happens  in the almost equal case. Many results are known in this case. For volume we have the following beautiful stability results for positive and nonnegative Ricci  curvatures \cite{Cheeger-Colding1997}.
\begin{theo}[Volume Stability, Cheeger-Colding, 1997]
There exists $\epsilon (n) >0$ such that \\
(i) if a complete Riemannian manifold $M^n$  has $\Ric_M \ge n-1$ and $\vol_M \ge (1 -\epsilon(n)) \vol(S^n)$, then $M^n$ is diffeomorphic to $S^n$; \\
(ii) if a complete Riemannian manifold $M^n$  has $\Ric_M \ge 0$ and for some $p \in M$, $\vol B(p,r) \ge (1 -\epsilon(n)) \, \omega_n r^n$ for all $r>0$, then $M^n$ is diffeomorphic to $\mathbb R^n$. \end{theo}
This was first proved by Perelman \cite{Perelman1994} with the weaker conclusion that $M^n$ is homeomorphic to $S^n$ (contractible resp.).

The analogous stability result is not true for diameter. In fact, there are
manifolds with $\Ric \ge n-1$ and diameter arbitrarily close to $\pi$
which are not homotopic to sphere \cite{Anderson1990-diameter, Otsu1991}. This should be contrasted with
 the sectional curvature case, where we have the beautiful Grove-Shiohama diameter sphere
theorem \cite{Grove-Shiohama1977}, that if
$M^n$ has sectional curvature $K_M \geq 1$ and $\diam_M >\pi/2$ then
$M$ is homeomorphic to $S^n$. Anderson showed that the stability  for the splitting theorem (Theorem~\ref{splitting}) does not hold either \cite{Anderson1992}.

By work of Cheng and Croke \cite{Cheng1975, Croke1982}, if $\Ric_M \ge
n-1$ then $\diam_M$ is close to $\pi$ if and only if $\lambda_1(M)$
is close to $n$. So the naive version of the stability for
$\lambda_1(M)$ does not hold either. However from the work of \cite{Colding1996-Large,
Cheeger-Colding1997, Petersen1999} we have the following modified
version.
\begin{theo}[Colding, Cheeger-Colding, Petersen]
There exists $\epsilon (n) >0$ such that if a complete Riemannian
manifold $M^n$  has $\Ric_M \ge n-1$, and radius $\ge \pi
-\epsilon(n)$ or  $\lambda_{n+1}(M) \le n+\epsilon (n)$, then $M^n$
is diffeomorphic to $S^n$.
\end{theo}
Here $\lambda_{n+1}(M)$ is the $(n+1)-$th eigenvalue of the Laplacian. The above condition is natural in the sense that for $S^n$
the radius is $\pi$ and the first eigenvalue is $n$ with multiplicity 
$n+1$. Extending Cheng and Croke's work Petersen showed that if $\Ric_M \ge n-1$ then the
radius is close to $\pi$ if and only if $\lambda_{n+1}(M)$ is close to
$n$.

The stability for the first Betti number,  conjectured by Gromov, was proved by Cheeger-Colding  in \cite{Cheeger-Colding1997}. Namely there exists $\epsilon (n) >0$ such that if a complete Riemannian
manifold $M^n$ has $\Ric_M (\diam_M)^2 \ge -\epsilon(n)$ and $b_1=n$,
then $M$ is diffeomorphic to $T^n$. The homeomorphic version was first proved in \cite{Colding1997}.

Although the direct stability for diameter does not hold, Cheeger-Colding's breakthrough work \cite{Cheeger-Colding1996} gives quantitative generalizations of  the diameter rigidity results, see Section~\ref{almost}.

 \sect{The fundamental groups}

In lower dimensions ($n \le 3$) a Ricci curvature lower bound has strong topological implications. R. Hamilton \cite{Hamilton1982} proved that compact manifolds $M^3$ with positive Ricci curvature are space forms.  Schoen-Yau \cite{Schoen-Yau1982} proved that any complete open manifold $M^3$ with positive Ricci curvature must be diffeomorphic to $\mathbb R^3$ using minimal surfaces. In general the strongest control is on the fundamental group.

The first  result is Myers' theorem \cite{Myers1941}.

\begin{theo}[Myers, 1941] \label{Myers} If $\Ric_M \geq H >0,$ then $\diam (M) \le \pi/\sqrt{H}$, and $\pi_1 (M)$ is finite.
\end{theo}

This is the only known topological obstruction to a compact manifold supports a metric with positive Ricci curvature other than topological obstructions shared by manifolds with positive scalar curvature. See Section~\ref{examples} for examples with positive Ricci curvature and
Rosenberg's article in this volume for a discussion of scalar curvature.

 We can still ask what one can say about the finite group.  Any finite group can be realized as the fundamental group of a
 compact manifold with positive Ricci curvature since any finite group is a subgroup of $SU(n)$ (for n sufficiently big) and $SU(n)$  has a metric with
 positive Ricci curvature (in fact Einstein).
 What can one say if the dimension $n$ is fixed? For example, is the order of the group modulo an abelian subgroup
 bounded by the dimension? See \cite{Wilking2000} for a partial result.

For a compact manifold $M$ with nonnegative Ricci curvature,  Cheeger-Gromoll's
splitting theorem (Theorem~\ref{splitting}) implies that  $\pi_1 (M)$ has an abelian subgroup
of finite index \cite{Cheeger-Gromoll1971}. Again it is open if one can bound the index by dimension.

For general nonnegative Ricci curvature manifolds, using covering and volume
comparison Milnor showed that \cite{Milnor1968}
 \begin{theo}[Milnor, 1968] If $M^n$ is complete with $\Ric_M \geq 0,$ then any
finitely generated subgroup of $\pi_1 (M)$ has polynomial growth
of degree $\leq n.$
\end{theo}

Combining this with the following result of Gromov \cite{Gromov1981},
we know that any finitely generated subgroup of $\pi_1 (M)$ of
manifolds with nonnegative Ricci curvature is almost nilpotent.
\begin{theo}[Gromov, 1981]  A finitely generated group $\Gamma$ has
polynomial growth iff $\Gamma$ is almost nilpotent, i.e. it contains a nilpotent subgroup of finite index.
\end{theo}

When $M^n$ has nonnegative Ricci curvature and Euclidean volume growth (i.e. $\vol B(p,r) \ge cr^n$ for some $c>0$),  using a heat kernel estimate Li showed that $\pi_1 (M)$ is finite \cite{Li1986}. Anderson also derived this using volume comparison \cite{MR1046624}. Using the splitting theorem of Cheeger and Gromoll \cite{Cheeger-Gromoll1971} (Theorem~\ref{splitting}) on the universal cover Sormani showed that a noncompact manifold with positive Ricci curvature has the loops-to-infinity property \cite{Sormani2001}. As a consequence she showed that a noncompact manifold with positive Ricci curvature is simply connected if it is simply connected at infinity. See \cite{Shen-Sormani2001, Wylie-thesis} for  more applications of the loops-to-infinity property.

 From the above one naturally wonders if all nilpotent groups occur as the fundamental group of a complete non-compact manifold with
 nonnegative Ricci curvature. Indeed, extending  the warping product constructions in
 \cite{Nabonnand1980, Berard-Bergery1986}, Wei showed \cite{Wei1988} that 
any finitely generated torsion free nilpotent group could occur as fundamental group of a manifold with
 positive Ricci curvature. Wilking \cite{Wilking2000} extended this to any finitely generated almost
nilpotent group. This gives a very good understanding of the fundamental group of a manifold with nonnegative Ricci curvature except the following
 long standing problem regarding the finiteness of generators \cite{Milnor1968}.
\begin{conj}[Milnor, 1968]
The fundamental group of a manifold with nonnegative Ricci curvature
is finitely generated.
\end{conj}

There is some very good progress in this direction. Using short
generators and a uniform cut lemma based
on the excess estimate of Abresch and Gromoll \cite{Abresch-Gromoll1990} (see (\ref{excess}) ) Sormani  \cite{Sormani2000} proved that
if $\Ric_M \ge0$ and $M^n$ has small linear
diameter growth, then $\pi_1 (M)$ is finitely generated.  More precisely the small linear growth condition
is:
 \[\limsup_{r \to \infty} \frac{\diam
\partial B(p,r)}{r} < s_n = \frac{n}{(n-1)3^n} \left(\frac{n-1}{n-2}
\right)^{n-1}. \]
The constant $s_n$ was improved in \cite{Xu-Wang-Yang2003}. Then in \cite{Wylie2006} Wylie proved  that in this case $\pi_1(M) = G(r)$ for $r$ big, where $G(r)$  is the image of $\pi_1(B(p, r))$ in $\pi_1(B(p, 2r))$. In an earlier paper \cite{Sormani-minvol}, Sormani proved that all manifolds with nonnegative
Ricci curvature and linear volume growth have sublinear diamter growth,
so manifolds with linear volume growth are covered by these results. Any open manifold
with nonnegative Ricci curvature has at least linear volume growth \cite{Yau1976}.

In a very different direction Wilking \cite{Wilking2000}, using algebraic methods,   showed that if  $\Ric_M \geq 0$ then $\pi_1 (M)$ is finitely generated iff any
abelian subgroup of $\pi_1 (M)$ is finitely generated,
effectively reducing the Milnor conjecture to the study of manifolds
with abelian fundamental groups.

The fundamental group and the first Betti number are very nicely related. So it is natural that Ricci lower bound also controls the first Betti number. For compact manifolds Gromov  \cite{Gromov1999} and Gallot \cite{Gallot1983} showed that if $M^n$ is a compact manifold with
\be \Ric_M \geq (n-1)H, \ \ \diam_M \leq D,  \label{ric-diam-bounds} \ee
then  there is a function $C(n, HD^2)$ such that $b_1 (M) \leq C(n, HD^2)$ and
$\displaystyle{\lim_{x \to 0^{-}}} C(n,x) = n$ and $C(n, x ) = 0 $
for $x > 0.$  In particular, if $HD^2$ is small, $b_1 (M) \leq n.$

The celebrated Betti number estimate of Gromov
\cite{Gromov1981-betti} shows that all higher Betti numbers can be
bounded by sectional curvature and diameter. This is not true for
Ricci curvature. Using semi-local surgery Sha-Yang constructed
metrics of positive Ricci curvature on the connected sum of $k$
copies of $S^2\times S^2$ for all $k \ge 1$ 
\cite{Sha-Yang1991}. Recently using Seifert bundles over orbifolds
with a K\"{a}hler Einstein metric Kollar showed that there are
Einstein metrics with positive Ricci curvature on the connected sums
of arbitrary number of copies of $S^2\times S^3$ \cite{Kollar}.

Kapovitch-Wilking \cite{Kapovitch-Wilking} recently announced a proof of 
the compact analog of Milnor's conjecture that the fundamental
group of a manifold satisfying (\ref{ric-diam-bounds}) has a
presentation with a universally bounded number of generators (as
conjectured by this author), and that a manifold which
 admits almost nonnegative Ricci curvature has a virtually nilpotent
fundamental group. The second result would greatly generalize Fukaya-Yamaguchi's work on almost nonnegative sectional curvature \cite{Fukaya-Yamaguchi1992}. See \cite{Wei1990, Wei1997} for earlier partial results.

When the volume is also bounded from below,  by using a clever covering argument  M. Anderson \cite{Anderson1990} showed that the number of the short homotopically nontrivial closed geodesics  can be controlled and for  the class of manifolds $M$ with $\Ric_M
\geq (n-1)H,$ $\vol_M \geq V$ and $\diam_M \leq D$ there are only
finitely many isomorphism types of $\pi_1 (M)$.
Again if the Ricci curvature is replaced by sectional curvature then much more can be said. Namely there are only finitely many homeomorphism types of the manifolds with sectional curvature and volume bounded from below and diameter bounded from above \cite{Grove-Petersen-Wu1990, Perelman-preprint}.  By \cite{Perelman1997} this is not true for Ricci curvature unless the dimension is 3 \cite{Zhu1993}.

Contrary to a Ricci curvature lower bound, a Ricci curvature upper bound does not have any topological constraint \cite{Lohkamp1994}.
\begin{theo}[Lohkamp, 1994]
If $n \ge 3$, any manifold, $M^n$, admits a complete metric with $\Ric_M <0$.
\end{theo}

An upper Ricci curvature bound does have geometric implications, e g. the isometry group of a compact manifold with negative Ricci curvature is finite. In the presence of a lower bound, an upper bound on Ricci curvature forces additional regularity of the metric, see Theorem~\ref{twoside} in Section~\ref{limit} by Anderson.
It's still unknown whether it will give additional topological control. For example, the following question is still open.
\begin{ques}
Does the class of manifolds $M^n$ with $|\Ric_M | \le H, \vol_M \ge V$ and $\diam_M \le D$ have finite many homotopy types?
\end{ques}
There are infinitely many homotopy types without the Ricci upper bound \cite{Perelman1997} .

 \sect{Gromov-Hausdorff convergence}  \label{Gromov-Hausdorff}

Gromov-Hausdorff convergence is very useful in studying manfolds
with a lower Ricci bound. The starting point is Gromov's
precompactness theorem. Let's first recall the Gromov-Hausdorff distance. See \cite[Chapter 3,5]{Gromov1999},\cite[Chapter 10]{Petersen-book}, \cite[Chapter 7]{Burago-Burago-Ivanov2001} for more background material on Gromov-Hausdorff convergence.

Given  a metric space $(X,d)$ and subsets $A,B \subset X$, the Hausdorff distance is
\[ d_H (A,B)=\inf\{\epsilon>0: B \subset T_\epsilon(A) \textrm{ and }
A \subset T_\epsilon(B) \},
\]
where $T_\epsilon(A)=\{x\in X: d (x,A) <\epsilon\}$.

\begin{defn}[Gromov, 1981] \label{GH-distance}
Given two compact metric spaces $X,Y$, the Gromov-Hausdorff distance  is
$ d_{GH}(X,Y) = \inf \left\{ d_H(X,Y): \right.$   all metrics on the disjoint union, $X \coprod Y$, which extend the metrics of $X$ and $\left.Y \right\}$.
\end{defn}

The Gromov-Hausdorff distance defines a metric on the collection of isometry classes of compact metric spaces. Thus,  there is the naturally associated notion of Gromov-Hausdorff convergence of
compact metric spaces.  While the Gromov-Hausdorff distance make sense for non-compact metric spaces,  the following looser definition of convergence is more appropriate. See also \cite[Defn 3.14]{Gromov1999}.  These two definitions are equivalent \cite[Appendix]{Sormani-Wei2004}.
\begin{defn}
We say that non-compact metric spaces $(X_i, x_i)$ converge in
the pointed Gromov-Hausdorff sense to $(Y,y)$ if
for any $r>0$, $B(x_i, r)$ converges to $B(y, r)$ in the pointed Gromov-Hausdorff
sense.
\end{defn}

Applying the relative volume comparison (\ref{relative-vol}) to  manifolds with lower Ricci bound,  we have
\begin{theo}[Gromov's  precompactness  theorem] \label{precompactness}
The class of closed manifolds $M^n$ with $\Ric_M \ge (n-1)H$ and $\diam_M \le D$ is precompact.
The class of pointed complete manifolds $M^n$ with $\Ric_M \ge (n-1)H$ is precompact.
\end{theo}

By the above, for an open manifold $M^n$ with $\Ric_M \ge 0$ any
sequence $\{ (M^n, x,r_i^{-2}g)\}$, with $r_i \ra \infty$,
subconverges in the pointed Gromov-Hausdorff topology to a length
space $M_\infty$. In general, $M_\infty$ is not unique
\cite{Perelman-cone}. Any such limit is called an asymptotic cone of
$M^n$, or a cone of $M^n$ at infinity .

Gromov-Hausdorff convergence defines a very weak topology. In general one
only knows that Gromov-Hausdorff limit of length spaces is a length
space and diameter is continuous under the Gromov-Hausdorff
convergence. When the limit is a smooth manifold with same dimension Colding showed the remarkable result that for manifolds
with lower Ricci curvature bound the volume also converges
\cite{Colding1997} which was conjectured by Anderson-Cheeger. See also \cite{Cheeger2001} for a proof using mod $2$ degree.
\begin{theo}[Volume Convergence, Colding, 1997]   \label{vol-conv} If $(M_i^n, x_i)$ has $\Ric_{M_i} \ge (n-1)H$ and converges in the pointed Gromov-Hausdorff sense to smooth Riemannian  manifold $(M^n,x)$, then for all $r >0$
\be
\lim_{i \ra \infty} \vol (B(x_i, r)) = \vol (B(x,r)).
\ee
\end{theo}

The volume convergence can be generalized  to the noncollapsed singular limit space (by replacing the Riemannian volume with the $n$-dimensional Hausdorff measure $\mathcal H^n$) \cite[Theorem 5.9]{Cheeger-Colding1997},  and to the collapsing case with smooth limit $M^k$ in terms of the $k$-dimensional Hausdorff content  \cite[Theorem 1.39]{Cheeger-Colding2000II}.

As an application of Theorem~\ref{vol-conv}, Colding \cite{Colding1997} derived the rigidity result that
if $M^n$ has $\Ric_{M} \ge 0$ and some $M_\infty$ is isometric to $\mathbb R^n$, then $M$ is isometric to $\mathbb R^n$.

We also have the following wonderful stability result \cite{Cheeger-Colding1997} which sharpens an earlier version in \cite{Colding1997}.
\begin{theo}[Cheeger-Colding, 1997]
For a closed Riemannian manifold $M^n$ there exists an $\epsilon(M) >0$ such that if $N^n$ is a $n$-manifold with $\Ric_N \ge -(n-1)$ and $d_{GH} (M,N) < \epsilon$ then $M$ and $N$ are diffeomorphic.
\end{theo}
Unlike the sectional curvature case, examples show that the result
does not hold if one allows $M$ to have singularities even
on the fundamental group level \cite[Remark (2)]{Otsu1991}. Also the
$\epsilon$ here must depend on $M$ \cite{Anderson1990-diameter}.

Cheeger-Colding also showed that the eigenvalues and eigenfunctions
of the Laplacian are continuous under measured Gromov-Hausdorff
convergence \cite{Cheeger-Colding2000III}. To state the result we
need a definition and some structure result on the limit space (see
Section~\ref{limit} for more structures). Let $X_i$ be a sequence of
metric spaces converging to $X_\infty$ and $\mu_i, \mu_\infty$ are
Radon measures on $X_i, X_\infty$.
\begin{defn}
We say $(X_i, \mu_i)$ converges in the measured Gromov-Hausdorff
sense to $(X_\infty, \mu_\infty)$ if for all sequences of continuous
functions $f_i:X_i \ra \mathbb R$ converging to $f_\infty:X_\infty
\ra \mathbb R$, we have
\be \int_{X_i} f_i d\mu_i \ra \int_{X_\infty} f_\infty d\mu_\infty. \ee
\end{defn}

 If $(M_\infty,p)$ is the pointed Gromov-Hausdorff limit of a sequence
of Riemannian manifolds $(M^n_i, p_i)$ with $\Ric_{M_i} \ge -(n-1)$,
then there is a natural collection of measures, $\mu$, on $M_\infty$ obtained by taking limits of the normalized Reimannian measures on $M^n_j$ for a suitable subsequence $M_j^n$ \cite{Fukaya1987}, \cite[Section 1]{Cheeger-Colding1997},
\be  \label{renormalized-vol}
\mu = \lim_{j \ra \infty} \underline{\vol}_j ( \cdot ) = \vol (\cdot )/\vol(B(p_j,1)).
\ee
  In particular, for all $z \in M_\infty$ and $0< r_1 \le r_2$, we have the renormalized limit measure $\mu$ satisfy the following comparison
  \be
  \frac{\mu (B(z,r_1))}{\mu (B(z, r_2))} \ge \frac{\vol_{n,-1}\left(B(r_1)\right)}{\vol_{n,-1}\left(B(r_2)\right) }. \ee
With this,  the extension of the segment inequality
(\ref{segment}) to the limit, the gradient estimate
(\ref{gradient2}), and Bochner's formula, one can define
a canonical self-adjoint Laplacian $\Delta_\infty$ on the limit
space $M_\infty$ by means of limits of the eigenfunctions and
eigenvalues for the sequence of the manifolds. In \cite{Cheeger1999,
Cheeger-Colding2000III}  an intrinsic construction of this operator
is also given on a more general metric measure spaces. Let
$\{\lambda_{1,i} \cdots,\}, \{\lambda_{1,\infty}, \cdots, \}$ denote
the eigenvalues for $\Delta_i, \Delta_\infty$ on $M_i, M_\infty$,
and $\phi_{j,i}, \phi_{j,\infty}$ the eigenfunctions of the jth
eigenvalues $\lambda_{j,i}, \lambda_{j,\infty}$. In
\cite{Cheeger-Colding2000III}  Cheeger-Colding in particular proved
the following theorem, establishing Fukaya's conjecture
\cite{Fukaya1987}.
\begin{theo}[Spectral Convergence, Cheeger-Colding, 2000] Let  $(M_i^n, p_i, \underline{\vol}_i )$ with $\Ric_{M_i} \ge -(n-1)$ converges to $(M_\infty,p,\mu)$ under measured Gromov-Hausdorff sense and $M_\infty$ is compact. Then for each $j$, $\lambda_{j,i} \ra \lambda_{j,\infty}$ and $\phi_{j,i} \ra  \phi_{j,\infty}$ uniformly as $i \ra \infty$.
\end{theo}

As a natural extension, in \cite{Ding2002} Ding proved that the heat kernel and Green's function also behave nicely  under the measured Gromov-Hausdorff convergence. The natural extension to the $p$-form Laplacian does not hold, however, there is still very nice work in this direction by John Lott, see \cite{Lott2002, Lott2004}.

 \sect{Almost rigidity and applications}  \label{almost}

Although the analogous stability results for maximal diameter in the case of
positive/nonnegative Ricci curvature do not hold, Cheeger-Colding's significant
work \cite{Cheeger-Colding1996} provides quantitative generalizations
of Cheng's maximal diameter theorem, Cheeger-Gromoll's splitting theorem (Theorem~\ref{splitting}), and the volume annulus
implies metric annulus theorem in terms of Gromov-Hausdroff
distance. These results have important applications in extending rigidity results  to
the limit space.

An important ingredient for these results is Abresch-Gromoll's excess
estimate \cite{Abresch-Gromoll1990}. For $y_1, y_2 \in M^n$, the excess
function $E$ with respect to $y_1,y_2$ is
 \be E_{y_1,y_2}(x) = d (y_1,x)+d(y_2,x) - d(y_1,y_2). \ee
 Clearly $E$ is Lipschitz with Lipschitz
constant $\le 2$.

 Let $\gamma$ be a minimal geodesic from $y_1$ to $y_2$, $s(x) = \min (d(y_1,x), d(y_2,x))$ and $h(x)=\min_{t} d(x, \gamma(t))$, the height from $x$ to a minimal geodesic $\gamma(t)$  connecting $y_1$ and $y_2$. By the triangle inequality $0 \le E(x) \le 2 h(x)$.
Applying the Laplacian comparison (Theorem~\ref{laplacian-comp}) to
$E(x)$ and with an elaborate (quantitative) use of the maximum
principle (Theorem~\ref{maximum-p}) Abresch-Gromoll showed that if
$\Ric_M \ge 0$ and $h(x) \le \frac{s(x)}{2}$, then
(\cite{Abresch-Gromoll1990}, see also \cite{Cheeger1991})
 \be E(x) \le 4 \left(\frac{h^{n}}{s}\right)^{\frac{1}{n-1}}.  \label{excess} \ee
This is the first distance estimate in terms of a lower Ricci curvature bound.

 The following version (not assuming $E(p) = 0$,
but without the sharp estimate) is from \cite[Theorem
9.1]{Cheeger2001}.

\begin{theo}[Excess Estimate, Abresch-Gromoll, 1990] \label{Abresch-Gromoll} If $M^n$ has $\Ric_M \ge -(n-1)\delta$, and for $p \in M$, $s(p) \ge L$ and $E(p) \le \epsilon$, then on $B(p, R)$, $E \le \Psi = \Psi(\delta,L^{-1},\epsilon |\,n,R)$, where $\Psi$ is a nonnegative constant such that for fixed $n$ and $R$ $\Psi$ goes to  zero as $\delta,\epsilon \ra 0$ and $L \ra \infty$.
\end{theo}

 This can be interpreted as a weak almost splitting theorem. Cheeger-Colding generalized this result tremendously by proving the following almost splitting theorem \cite{Cheeger-Colding1996}, see also \cite{Cheeger2001}.
\begin{theo}[Almost Splitting, Cheeger-Colding, 1996] \label{almost-splitting} With the same assumptions as Theorem~\ref{Abresch-Gromoll}, there is a
length space $X$ such that for some ball, $B((0,x), \frac{1}{4}R) \subset \mathbb R \times X$, with the product metric, we have
\[
d_{GH} \left( B(p, \frac{1}{4}R), B((0,x), \frac{1}{4}R) \right) \le \Psi.
\]
\end{theo}
Note that $X$ here may not be smooth, and the Hausdorff dimension could be smaller
than $n-1$. Examples also show that the ball $B(p, \frac{1}{4}R)$
may not have the topology of a product, no matter how small $\delta,
\epsilon$, and $L^{-1}$  are \cite{Anderson1992,
Menguy2000-noncollapsing}.

The proof is quite involved. Using the Laplacian comparison,
the maximum principle, and Theorem~\ref{Abresch-Gromoll}  one shows that
the distance function $b_i = d(x,y_i)-d(p,y_i)$  associated to $p$ and
$y_i$ is uniformly close to $\mathbf b_i$, the harmonic function with same
values on $\partial B(p,R)$. From this, together with the lower bound for the smallest eigenvalue of the Dirichlet problem on $B(p,R)$ (see Theorem~\ref{thm-segment}) one shows that $\nabla b_i, \nabla \mathbf b_i$ are close in the $L_2$ sense. In particular $\nabla \mathbf b_i$ is close to $1$ in the $L_2$ sense. Then  applying the Bochner formula to $\mathbf b_i$ multiplied with a cut-off function with bounded Laplacian one shows that
  $| \Hess \mathbf b_i|$ is small in the $L_2$ sense in a
smaller ball. Finally, in the most significant step, by using the segment
inequality (\ref{segment}), the gradient estimate (\ref{gradient}) and the information established above one derives a quantitative version of the Pythagorean
theorem, showing that the ball is close in the Gromov-Hausdorff sense to a ball in some product space; see \cite{Cheeger-Colding1996, Cheeger2001}.

An immediate application of the almost splitting theorem is the extension of the splitting theorem to the limit space.
\begin{theo}[Cheeger-Colding, 1996]  If $M_i^n$ has $\Ric_{M_i} \ge -(n-1)\delta_i$ with $\delta_i \ra 0$ as $i\ra \infty$, converges to $Y$ in the pointed Gromov-Hausdorff sense, and $Y$ contains a line, then $Y$ is isometric to $\mathbb R \times X$ for some length space $X$.
\end{theo}

Similarly, one has almost rigidity in the presence of finite diameter (with simpler a proof)
\cite[Theorem 5.12]{Cheeger-Colding1996}. As a special consequence, we have that if $M_i^n$ has $\Ric_{M_i} \ge (n-1)$,  $\diam_{M_i} \ra \pi$ as $i\ra \infty$, and converges to $Y$ in the Gromov-Hausdorff sense, then $Y$ is isometric to the spherical metric suspension of some length space $X$ with $\diam (X) \le \pi$. This is a kind of stability for diameter. 

Along the same lines (with more complicated technical details) Cheeger and Colding \cite{Cheeger-Colding1996} have an almost rigidity version for the volume annulus implies metric annulus theorem (see Section 5). As a very nice application to the
asymptotic cone, they showed that if $M^n$ has $\Ric_{M} \ge 0$  and has Euclidean volume growth, then every asymptotic cone of $M$ is a metric cone.

\sect{The structure of limit spaces}  \label{limit}

As we have seen, understanding the structure of the limit space of manifolds with
lower Ricci curvature bound often helps in understanding the
structure of the sequence. Cheeger-Colding made significant progress
in understand the regularity and geometric structure of the limit
spaces \cite{Cheeger-Colding1997, Cheeger-Colding2000II,
Cheeger-Colding2000III}. On the other hand Menguy constructed examples showing that
the limit space could have infinite topology in an arbitrarily small neighborhood \cite{Menguy2000-noncollapsing}. In \cite{Sormani-Wei2001,Sormani-Wei2004}
Sormani-Wei showed that the limit space has a universal cover.

Let $(Y^m,y)$ (Hausdorff dimension $m$) be the pointed Gromov-Hausdorff limit of a sequence of Riemannian manifolds $(M^n_i, p_i)$ with $\Ric_{M_i} \ge -(n-1)$. Then $m \le n$ and $Y^m$ is locally compact. Moreover Cheeger-Colding \cite{Cheeger-Colding1997} showed that if $m = \dim Y < n$, then $m \le n-1$.

The basic notion for studying the infinitesimal structure of the limit space $Y$ is that of a tangent cone.
\begin{defn}
A tangent cone, $Y_y$, at $y\in (Y^m,d)$ is the pointed Gromov-Hausdorff limit of a sequence of the rescaled spaces $(Y^m,r_id,y)$, where $r_i \ra \infty$ as $i \ra \infty$.
\end{defn}

By Gromov's precompactness theorem (Theorem~\ref{precompactness}),
every such sequence has a converging subsequence. So tangent cones
exist for all $y \in Y^m$, but might depend on the choice of
convergent sequence. Clearly if $M^n$ is a Riemannian manifold, then
the tangent cone at any point is isometric to $\mathbb R^n$.
Motivated by this one defines \cite{Cheeger-Colding1997}

\begin{defn} A point, $y\in Y$, is called  $k$-regular  if  for some $k$, every  tangent
cone at $y$ is isometric to $\mathbb R^k$. Let $\mathcal R_k$ denote the set of $k$-regular  points and $\mathcal R = \cup_k \mathcal R_k$, the regular set. The  singular set, $Y\setminus \mathcal R$, is denoted $\mathcal S$.
 \end{defn}

 Let $\mu$ be a renormalized limit measure on $Y$ as in (\ref{renormalized-vol}). Cheeger-Colding showed that the regular points have full measure \cite{Cheeger-Colding1997}.

\begin{theo}[Cheeger-Colding, 1997]   For any renormalized limit measure $\mu$, $\mu(\mathcal S) =0$, in particular, the regular points are dense.
\end{theo}

Furthermore, up to a set of measure zero, $Y$ is a countable union of
sets, each of which is bi-Lipschitz equivalent to a subset of
Euclidean space \cite{Cheeger-Colding2000III}.
\begin{defn}
A metric measure space, $(X, \mu)$, is called $\mu$-rectifiable if
$0<\mu(X) <\infty$, and there exists $N < \infty$ and a countable
collection of subsets, $A_j$, with $\mu(X\setminus \cup_jA_j) = 0$,
such that each $A_j$ is bi-Lipschitz equivalent to a subset of
$\mathbb R^{l(j)}$, for some $1 \le l(j) \le N$ and in addtion, on
the sets $A_j$, the measures $\mu$ and and the Hausdorff measure
$\mathcal H^{l(j)}$ are mutually absolutely continous.
\end{defn}

\begin{theo}[Cheeger-Colding, 2000]
Bounded subsets of $Y$ are $\mu$-rectifiable with respect to any
renormalized limit measure $\mu$.
\end{theo}

At the singular points, the structure could be very complicated. Following a related earlier construction of Perelman \cite{Perelman1997}, Menguy constructed 4-dimensional examples of (noncollapsed) limit spaces with, $\Ric_{M_i^n} >1$, for which there exists point so that any neighborhood of the point  has infinite second Betti number \cite{Menguy2000-noncollapsing}. See \cite{Cheeger-Colding1997, MR1781925, Menguy2001} for examples of collapsed limit space with interesting properties.

Although we have very good regularity results, not much topological structure is known for the limit spaces in general. E.g., is $Y$ locally simply connected? Although this is unknown, using the renormalized limit measure and the existence of regular points, together with $\delta$-covers, Sormani-Wei \cite{Sormani-Wei2001,Sormani-Wei2004} showed that
the universal cover of $Y$ exists. Moreover when $Y$ is compact, the fundamental group of $M_i$ has a surjective homomorphism  onto the group of deck
transforms of $Y$ for all $i$ sufficiently large.

When the sequence has the additional assumption  that
\be
\vol (B(p_i,1)) \ge v >0,
\ee
the limit space $Y$ is called noncollapsed. This is equivalent to $m=n$. In this case, more structure is known.
\begin{defn}
Given $\epsilon > 0$, the $\epsilon$-regular set, $\mathcal R_\epsilon$, consists of those points $y$ such that for all sufficiently small $r$,
\[
d_{GH} (B(y,r), B(0,r)) \le \epsilon r, \]
where $0 \in \mathbb R^n$.
\end{defn}
Clearly $\mathcal R = \cap_\epsilon \mathcal R_\epsilon$.  Let $\stackrel{\circ}{\mathcal R}_\epsilon$ denote the interior of $\mathcal R_\epsilon$.
\begin{theo}[Cheeger-Colding 1997, 2000] There exists $\epsilon (n) > 0$ such that if $Y$ is a noncollapsed limit space of the sequence $M_i^n$ with $\Ric_{M_i} \ge -(n-1)$, then for $0 < \epsilon < \epsilon(n)$, the set $\stackrel{\circ}{\mathcal R}_\epsilon$ is $\alpha(\epsilon)$-bi-H\"older equivalent to a smooth connected Riemannian manifold, where $\alpha(\epsilon) \ra 1$ as $\epsilon \ra 0$. Moreover,
\be
\dim(Y\setminus \stackrel{\circ}{\mathcal R}_\epsilon ) \le n-2.
\ee
In addition, for all $y \in Y$, every tangent cone $Y_y$ at $y$ is a metric cone and the isometry group of $Y$ is a Lie group.
\end{theo}

This is proved in \cite{Cheeger-Colding1997, Cheeger-Colding2000II}.

If, in addition, Ricci curvature is bounded from two sides, we have stronger regularity \cite{Anderson-MR1074481}.
\begin{theo}[Anderson, 1990]  \label{twoside} There exists $\epsilon (n) > 0$ such that if $Y$ is a noncollapsed limit space of the sequence $M_i^n$ with $| \Ric_{M_i} | \le n-1$, then for $0 < \epsilon < \epsilon(n)$, $\mathcal R_\epsilon = \mathcal R$. In particular the singular set is closed. Moreover, $\mathcal R$ is a $C^{1,\alpha}$ Riemannian manifold, for all $\alpha <1$. If the metrics on $M_i^n$ are Einstein, $\Ric_{M^n_i}  = (n-1)Hg_i$, then the metric on $\mathcal R$ is actually $C^\infty$.
\end{theo}

Many more regularity results are obtained when the sequence is Einstein, K\"ahler, has special holonomy, or has bounded $L^p$-norm of the full curvature tensor, see \cite{Anderson-Cheeger1991, Cheeger2003, Cheeger-Colding-Tian2002, Cheeger-Tian2005}, especially \cite{Cheeger2002} which gives an excellent survey in this direction. See the recent work \cite{Cheeger-Tian2006} for Einstein 4-manifolds with possible collapsing.

\sect{Examples of manifolds with nonnegative Ricci curvature} \label{examples}

Many examples of manifolds with nonnegative Ricci curvature have been constructed, which contribute greatly to the study of manifolds with lower Ricci curvature bound. We only discuss the examples related to the basic methods here, therefore many specific examples are unfortunately omitted (some are mentioned in the previous sections). There are mainly three methods: fiber bundle construction, special surgery,  and group quotient, all combined with warped products. These method are  also very useful in constructing Einstein manifolds. A large class of Einstein manifolds is also provided by Yau's solution of Calabi conjecture.

Note that if two compact Riemannian manifolds $M^m, N^n (n,m \ge 2)$ have  positive Ricci curvature, then their product has positive Ricci curvature, which is not true for sectional curvature but only needs one factor positive for scalar curvature. Therefore it is natural to look at the fiber bundle case.  Using Riemannian submersions with totally geodesic fibers J. C. Nash \cite{Nash1979}, W. A. Poor \cite{Poor1975}, and Berard-Bergery \cite{Berard-Bergery1978} showed that the compact total space of a fiber bundle admits a metric of positive Ricci curvature if the base and the fiber admit metrics with positive Ricci curvature and if the structure group acts by isometries. Furthermore, any vector bundle of rank $\ge 2$ over a compact manifold with $\Ric >0$ carries a complete metric with positive Ricci curvature. In \cite{Gilkey-Park-Tuschmann1998} Gilkey-Park-Tuschmann showed that a principal bundle $P$ over a compact manifold with $\Ric >0$ and compact connected structure group $G$  admits a $G$ invariant metric with positive Ricci curvature if and only if $\pi_1(P)$ is finite. Unlike the product case, the corresponding statements for $\Ric \ge 0$ are not true in all these cases, e.g. the nilmanifold $S^1 \ra N^3\ra T^2$ does not admit a metric with $\Ric \ge 0$. On the other hand  Belegradek-Wei \cite{Belegradek-Wei2004}  showed that it is true in the stable sense. Namely, if $E$ is the total space of a bundle over a compact base with $\Ric \ge 0$, and either a compact $\Ric \ge 0$  fiber or  vector space as fibers, with compact structure group acting by isometry, then $E\times \mathbb R^p$ admits a complete metric with positive Ricci curvature for all sufficiently large $p$. See \cite{Wraith} for an estimate of $p$.

 Surgery constructions are very successful in constructing manifolds with positive scalar curvature, see Rothenberg's article in this volume.  Sha-Yang  \cite{Sha-Yang1989, Sha-Yang1991}  showed that this is also a useful method  for constructing manifolds with positive Ricci curvature in special cases. In particular they showed that if $M^{m+1}$ has a complete metric with $\Ric >0$, and $n, m \ge 2$, then $S^{n-1} \times \left( M^{m+1} \setminus \coprod_{i=0}^k D_i^{m+1}\right) \bigcup_{Id} D^n \times \coprod_{i=0}^k S_i^{m}$, which is diffeomorphic to $\left( S^{n-1} \times M^{m+1}\right)  \# \left( \#_{i=1}^{k} S^n \times S^m \right)$, carries a complete metric with $\Ric >0$ for all $k$, showing that the total Betti number of a  compact Riemannian $n$-manifold ($n \ge 4$) with positive Ricci curvature could be arbitrarily large. See also \cite{Anderson1992}, and \cite{Wraith1998} when the gluing map is not the identity.

Note that a compact homogeneous space admits an invariant metric with positive Ricci curvature if and only if the fundamental group is finite \cite[Proposition 3.4]{Nash1979}. This is extended greatly by Grove-Ziller \cite{Grove-Ziller2002} showing that
any cohomogeneity one manifold $M$ admits a complete invariant metric
with nonnegative Ricci curvature and if $M$ is compact then it has positive Ricci curvature if
and only if its fundamental group is finite (see also \cite{Schwachhofer-Tuschmann2004}).   Therefore, the fundamental group is the only obstruction to a compact manifold admitting a positive Ricci curvature metric when there is enough symmetry.
It remains open what the obstructions  are to positive Ricci curvature besides the restriction on the fundamental group and those coming from positive scalar curvature (such as the $\hat A$-genus).

Of course, another interesting class of examples are given by
Einstein manifolds. For these,  besides the ``bible"
on Einstein manifolds \cite{Besse1987}, one can refer to the survey
book \cite{Einstein-survey-book} for the development after
\cite{Besse1987}, and the recent articles
\cite{Boyer-Galicki-preprint, Bohm-Wang-Ziller2004} for Sasakian
Einstein metrics and  compact homogenous Einstein manifolds.

\def\cprime{$'$}

Department of Mathematics,

University of California,

Santa Barbara, CA 93106

wei@math.ucsb.edu
\end{document}